 \documentclass[10pt,reqno]{amsart}
  \usepackage{geometry}
\usepackage{amsfonts}
\usepackage{amssymb}
\usepackage{amscd}
\usepackage[dvips]{graphicx}
\usepackage{xcolor}

\usepackage[all]{xy}

\title{Differential graded Brauer groups}

\date{March 31, 2023}

\author{Alexander Zimmermann}
\address{
Universit\'e de Picardie,
D\'epartement de Math\'ematiques et LAMFA (UMR 7352 du CNRS),
33 rue St Leu,
F-80039 Amiens Cedex 1,
France}
\email{alexander.zimmermann@u-picardie.fr}


\newtheorem{Lemma1}{{Lemma}}[section]
\newtheorem{Def1}[Lemma1]{{Definition}}
\newtheorem{Prop1}[Lemma1]{{Proposition}}
\newtheorem{Claim1}[Lemma1]{{Claim}}
\newtheorem{Rem1}[Lemma1]{{Remark}}
\newtheorem{Cor1}[Lemma1]{{Corollary}}
\newtheorem{Ex1}[Lemma1]{{Example}}
\newtheorem{Qu1}[Lemma1]{{Question}}
\newtheorem{Theo1}[Lemma1]{{Theorem}}
\newtheorem*{Theo2}{{Theorem}}

\newenvironment{Lemma}{\begin{Lemma1}}{\end{Lemma1}}
\newenvironment{Def}{\begin{Def1}\em}{\end{Def1}}
\newenvironment{Prop}{\begin{Prop1}}{\end{Prop1}}
\newenvironment{Rem}{\begin{Rem1}\rm}{\end{Rem1}}
\newenvironment{Theorem}{\begin{Theo1}}{\end{Theo1}}

\newenvironment{Example}{\begin{Ex1}\em}{\end{Ex1}}

\newcommand{\lra}{\longrightarrow}

\newcommand{\ra}{\rightarrow}
\newcommand{\sdp}{\times\kern-.2em\vrule height1.1ex depth-.05ex}
\newcommand{\epi}{\lra \kern-.8em\ra}

\newcommand{\Q}{{\mathbb Q}}
\newcommand{\N}{{\mathbb N}}

\newcommand{\Z}{{\mathbb Z}}

\newcommand{\End}{\textup{End}}
\newcommand{\Hom}{\textup{Hom}}
\newcommand{\GBr}{\textup{GBr}}
\newcommand{\dgBr}{\textup{dgBr}}
\newcommand{\Br}{\textup{Br}}
\newcommand{\Mat}{\textup{Mat}}
\newcommand{\red}{}

\begin{document}

\begin{abstract}
We consider central simple $K$-algebras which happen to be
differential graded $K$-algebras. Two such algebras $A$ and $B$
are considered equivalent if there are bounded complexes of finite dimensional
$K$-vector spaces $C_A$ and $C_B$ such that the differential graded algebras
$A\otimes_K \End_K^\bullet(C_A)$ and $B\otimes_K \End_K^\bullet(C_B)$ are isomorphic.
Equivalence classes form an abelian group, which we call the dg Brauer group.
We prove that this group is isomorphic to the ordinary Brauer group of the field $K$.
\end{abstract}

\maketitle

\section*{Introduction}

Brauer groups proved to be an important invariant of fields $K$. They are defined as
equivalence classes of central simple $K$ algebras, where two such algebras are called equivalent
if a matrix algebra of the one is isomorphic to the matrix algebra (of possibly different size)
of the other. These then form an abelian group under tensor product over $K$. Most
interestingly, for algebraic number fields, is the link to the Brauer groups of the
completions at the primes of the field. The Brauer group then embeds into the
product of the Brauer groups over all completions, with cokernel being isomorphic to $\Q/\Z$.

The original definition was generalised further to graded  central simple $K$
algebras (originally by C.T.C. Wall for a grading over the cyclic group of
order $2$, motivated by studies on quadratic forms), to Hopf algebras, and
culminating in the maybe the most far reaching generalisation to braided monoidal categories,
which was given by
van Oystaen and Zhang \cite{vanOystaen}.
Most recently, \cite{Gupta} developed a theory of Brauer groups for
(non graded!) differential central simple algebras in the context of
differential Galois theory. However, this theory is very different since
skew field never contain nilpotent elements, and hence only trivial gradings on
finite dimensional skew fields are possible. Therefore, skew fields
only bear a trivial differential graded structure.

In \cite{dgorders} we studied differential graded orders in differential
graded algebras which are semisimple  as algebras. In particular we
studied the local-global behaviour and defined a theory of id\`eles of
theses structures. It seems to be natural then to consider Brauer groups
for differential graded algebras which are central simple as algebras.
This is what we propose to do in the present note.

There are two
possible versions. First we may consider differential graded $K$-algebras which are
central simple as algebras, which leads to what we call
$\textup{dgBr}(K)$, the differential graded Brauer group.

Another option is the construction given by algebras which are finite
dimensional differential graded
$K$-algebras and whose category of differential graded modules is semisimple.
Aldrich and Garcia-Rozas  \cite{Tempest-Garcia-Rochas} proved a
structure theorem for these.
They are formed by differential graded algebras $(A,d)$, which are acyclic
as complexes, and such that $\ker(d)$ is central $\Z$-graded simple.
Moreover, in this case for any $z\in d^{-1}(1)$ we
have $A=\ker(d)\oplus z\ker(d)$. This concept does not give a Brauer group
since the tensor product of two such algebras will in general not be
simple in the sense of \cite{Tempest-Garcia-Rochas}, and actually not
even semisimple. We shall provide a
counterexample.

Our main result is the proof that the forgetful functor induces a
group isomorphism
$\textup{dgBr}(K)\simeq\textup{Br}(K)$.

The paper is organized as follows. In Section~\ref{generalities} we recall the
definitions and notations which we use for differential graded algebras.
Section~\ref{Brauergroupsdef} then recalls the definition of the ordinary Brauer group
of central simple algebras, the definition of Brauer groups in the graded
sense, and gives our main definition, namely the Brauer group of differential
graded algebras and shows first properties. We also provide an example why
the theory of semisimplicity in the category of differential graded algebras given
by \cite{Tempest-Garcia-Rochas} is not well suited for our purposes.
Section~\ref{mainresultsection} then states and proves our main result.

\section{Foundations of dg-algebras and dg-modules}

\label{generalities}

\subsection{Generalities}

We recall some definitions concerning differential graded algebras and their
differential graded modules, also to fix the notations.

Let $R$ be a commutative ring. Recall from Cartan \cite{Cartandg}, Keller \cite{Kellerdgandtilting}, and \cite{Stacks} that
\begin{enumerate}
\item
a differential graded $R$-algebra (or dg-algebra for short) is a $\Z$-graded
$R$-algebra $A=\bigoplus_{n\in\Z}A_n$ together with a graded $R$-linear
endomorphism $d$ satisfying $d(A_n)\subseteq A_{n+1}$ and $d\circ d=0$, and such that $$d(ab)=d(a)b+(-1)^{|a|}ad(b)$$
for all homogeneous elements  $a,b\in A$. Here $|a|$ denotes the degree of $a$. A homomorphism of differential graded algebras
$f:(A,d_A)\ra (B,d_B)$ is a degree $0$ homogeneous $R$-algebra map $f$ such that
$f\circ d_A=d_B\circ f$.
\item
for differential graded $K$-algebras $(A,d_A)$ and $(B,d_B)$, also
$A\otimes_KB$ is a differential graded algebra. The differential is defined by 
$$d_{A\otimes_KB}(a\otimes b)=d_A(a)\otimes b+(-1)^{|a|}a\otimes d_B(b)$$
for homogeneous elements $a$ and $b$. 
\item
If $(A,d)$ is a differential graded algebra, then $(A^{op},d^{op})$ is a differential
graded algebra (cf e.g. \cite[Definition 11.1]{Stacks})
with $x\cdot_{op}y:=(-1)^{|x|\cdot|y|}yx$ for any homogeneous elements $x,y\in A$, and $d^{op}(x)=d(x)$.
We hence write $d^{op}=d$.
\item
A differential graded right $A$-module (or dg-module for short) is then a $\Z$-graded $R$-module $M$ with graded $R$-linear
endomorphism $d_M$ of square $0$ and of degree $1$, such that $$d_M(ma)=d_M(m)a+(-1)^{|m|}md(a)$$
for all homogeneous elements $a\in A$ and $m\in M$.
A differential graded $(A,d)$-left module is a differential graded $(A^{op},d)$-right module.
\item
Let $(A,d_A)$ be a differential graded $R$-algebra and let $(M,\delta_M)$ and $(N,\delta_N)$ be
differential graded $(A,d_A)$-modules. Then a homomorphism of differential
graded modules is an $R$-linear homogeneous map $f:M\rightarrow N$ of degree $0$
with $f\circ \delta_M=\delta_N\circ f$, with
$f(am)=af(m)$ for all $a\in A$ and $m\in M$.
\item \label{lastitem} Let $(M,d_M)$ and $(N,d_N)$ be
differential graded $(A,d)$-modules.
The homomorphism complex $\Hom_A^\bullet(M,N)$ is the $\Z$-graded $R$-module given by
$$(\Hom_A^\bullet(M,N))_n:=\{f:M\ra N\;|\;f\in \Hom_R(M,N)\;\mbox{ and }f(M_k)\subseteq N_{k+n}\}.$$
The elements $f$ of $\Hom_A^\bullet(M,N)$  are not asked to be compatible with the differentials in any way.
Let $d_{\Hom}: \Hom_A^\bullet(M,N)\ra \Hom_A^\bullet(M,N)$ given by
$$d_{\Hom}(f):=d_N\circ f-(-1)^{|f|}f\circ d_M.$$
Then $d_{\Hom}^2(f)=0$, as is easily verified (cf e.g.~\cite{dgorders}).
Hence the pair $(\Hom_A^\bullet(M,N),d_{\Hom})$ is a complex of $R$-modules.
Moreover, in case $M=N$ we get that $(\Hom_A^\bullet(M,M),d_{\Hom})$ is a differential graded algebra.
\end{enumerate}

\section{Definition of differential graded Brauer groups}

\label{Brauergroupsdef}

We recall the Brauer group of a field.

\begin{Def}
Let $K$ be a field. Two finite dimensional central simple $K$-algebras $A$ and $B$ are equivalent if
there are positive integers $m,n\in\N$ such that $A\otimes_K\Mat_n(K)\simeq B\otimes_K\Mat_m(K)$.
The {\em Brauer group $\Br(K)$} is the group with elements being the equivalence classes of finite
dimensional central simple $K$-algebras and group law induced by $-\otimes_K-$.
\end{Def}

Recall that this is indeed a group. It is a set since we consider equivalence classes. 
The law is clearly well-defined and associative. The neutral element
is the equivalence class of $K$. Further,
for any finite dimensional central simple $K$-algebra $A$ we have
a classical result (cf e.g. \cite[Theorem 4.1.3]{Herstein}), which shows that
$A\otimes_KA^{op}\simeq \Mat_{\dim_KA}(K)$. Hence the inverse of
the equivalence class of $A$ is the equivalence class of $A^{op}$.

\medskip

For the so-called $G$-graded Brauer group $\GBr(K)$, for the (possibly infinite) cyclic group $G$,
an analogous construction is used. Note however that the graded Brauer group is
slightly different from the one in Definition~\ref{dgbrauerdef}.

A $G$-graded $K$-algebra $A$ is called graded central simple if $A$ does not have
non trivial $G$-graded twosided ideals and $K=Z(A)$. If $A$ and $B$ are two
graded central simple $K$-algebras, then $A\otimes_KB$ is again a
graded central simple $K$-algebra. Note that here the tensor product is the graded tensor product
defined as follows.
As a $K$-module, this is the usual tensor product. However the multiplication law
is given by $(a\otimes b)\cdot (c\otimes d)=(-1)^{|b||c|}(ac\otimes bd)$
for homogeneous elements $a,b,c,d$.

Again, it can be shown that $A\otimes_KA^{op}\simeq \End_K(A)$ as a graded algebra. As in the
ungraded case one defines an abelian group, the $G$-graded Brauer group (cf e.g. Turbow~\cite{Turbow}).

\begin{Def}
Let $K$ be a field and let $G$ be a cyclic group.
Two central $G$-graded simple $K$-algebras
$A$ and $B$ are equivalent if there are $G$-graded $K$-modules $V$ and $W$ such that
$A\otimes_K\End_K(V)\simeq B\otimes_K\End_K(W)$. Then the {\em $G$-graded Brauer group
$\GBr(K)$} is the group of equivalence classes of central $G$-graded simple $K$-algebras.
\end{Def}

For the $G$-graded Brauer group, one considers algebras
which are simple as $G$-graded algebras. Algebras which are simple as $G$-graded
algebras need not be simple as algebras. However, simple algebras, which happen
to be $G$-graded, are of course graded simple. Obviously simple algebras
are graded simple with the trivial grading. The graded
Brauer group is well-studied, and in case $G$ is or order $2$, this group is called the
Brauer-Wall group, after C.T.C. Wall's work~\cite{CTCWall}.

\subsection{General properties on differential graded algebras}

Our intention is to define a Brauer group for differential graded algebras.

We have (at least) two concepts for what we should call a simple differential graded algebra.
The somehow naive version, but underpinned by the success of the concept of
differential graded orders in \cite{dgorders} consists in considering central simple $K$-algebras
$A$, which are in addition differential graded. Our objects then would be such algebras $(A,d)$.

A second more categorical concept would be to consider indecomposable differential graded
algebras whose category of differential graded modules is semisimple.

Aldrich and Garcia-Rochas showed in \cite{Tempest-Garcia-Rochas} that
a differential graded algebra $(A,d)$ has a semisimple category of 
differential graded modules if and only if
$(A,d)$ is bounded and acyclic, and moreover $\ker(d)$ is semisimple as an ordinary graded
algebra. Further, in this case, for $z\in A$ with $d(z)=1$, we have $A=\ker(d)\oplus z\cdot\ker(d)$.
Observe that for any homogeneous element $n\in\ker(d)$ we have
$$d(zn)=d(z)n+(-1)^{|z|}z\cdot d(n)=1\cdot n+0=n.$$
Note further that we may use the same element $z$ also for the opposite algebra.

\begin{Example} \label{examplesimpledg}
For a field $K$ 
consider $A=K[X]/X^2$ where $K$ is in
degree $0$ and $X$ is an element of degree $-1$.
Then there is a differential $d$ on $K[X]/X^2$ given by $d(1)=0$ and $d(X)=1$.
Note that $(A,d)$ is acyclic and $\ker(d)=K$ is simple. Hence $(A,d)$ is semisimple in the
sense of \cite{Tempest-Garcia-Rochas}, but $A$ is not semisimple as an algebra.

On the other hand, the field $K$ in degree $0$, and differential $0$ is a
simple algebra, which is differential graded by the trivial dg-structure. However, it is not
simple in the sense of \cite{Tempest-Garcia-Rochas}.
\end{Example}

We hence first
need an analogue for what provides the inverse of a central simple algebra in the classical case.

\begin{Prop} \label{inverseofA}
Let $K$ be a field.
Let $(A,d)$ be a differential graded $K$-algebra which is central simple as an algebra.
Then
$(A,d)\otimes_K(A^{op},d_{A}^{op})\simeq (\End^\bullet_K(A),d_\Hom)$, and hence
$(A,d)\otimes_K(A^{op},d_{A}^{op})\simeq(\Mat_{\dim_K(A)}(K),d_M)$ for some
grading and differential $d_M$
on the matrix algebra.
\end{Prop}

\begin{proof} As $(A,d)$ is a (the regular) $(A,d)$-module, by the preliminary remarks,
$(A,d)$ is a right $(A^{op},d)$-module as well.
Hence, considering $(A,d)$ as a bounded complex $(C,d)$ of $K$-modules,
we consider left multiplication
$$\lambda:(A,d)\lra (\End_K^\bullet(C),d_\Hom)$$
given by $\lambda_a(x)=ax$ for any $a\in (A,d)$ and $x\in (C,d)$
as well as right multiplication given by
$$\rho:(A^{op},d)\lra (\End_K^\bullet(C),d_\Hom)$$
given by $\rho_a(x)=a\cdot_{op}x$ for any $a\in (A^{op},d)$ and $x\in (C,d)$.
These are clearly algebra homomorphisms. We need to see that these
are dg-homomorphisms between dg-algebras.

We first need to show $\lambda_{d(a)}=d_\Hom(\lambda_a)$ for all homogeneous $a\in A$. This
then translates into
\begin{align*}
d_\Hom(\lambda_a)(b)=&(d\circ\lambda_a-(-1)^{|a|}\lambda_a\circ d)(b)\\
=&d(ab)-(-1)^{|a|}a\cdot d(b)\\
=&d(a)\cdot b\\
=&\lambda_{d(a)}(b)
\end{align*}
for all $b\in A$.

Then, we need to show
$\rho_{d(a)}=d_\Hom(\rho_a)$ for all homogeneous $a\in A^{op}$.
This then becomes
\begin{align*}
d_\Hom(\rho_a)(b)=&(d\circ\rho_a-(-1)^{|a|}\rho_a\circ d)(b)\\
=&d(\rho_a(b))-(-1)^{|a|}\rho_a(d(b))\\
=&d(a\cdot_{op}b)-(-1)^{|a|}a\cdot_{op}d(b)\\
=&(-1)^{|b||a|}d(ba)-(-1)^{|a|}\cdot (-1)^{|a|\cdot(|b|+1)}d(b)a\\
=&(-1)^{|b||a|}(d(ba)-(-1)^{|a|+|a|}d(b)a)\\
=&(-1)^{|b||a|}(d(ba)-d(b)a)\\
=&(-1)^{|b||a|}\cdot(-1)^{|b|}bd(a)\\
=&(-1)^{|b|\cdot(|a|+1)}bd(a)\\
=&d(a)\cdot_{op}b\\
=&\rho_{d(a)}(b)
\end{align*}

Clearly $\rho$ is injective, as well as $\lambda$. Hence
$$A\simeq \lambda(A)\subseteq (\End_K^\bullet(C),d_\Hom)$$
and $$A^{op}\simeq \rho(A^{op})\subseteq (\End_K^\bullet(C),d_\Hom).$$
Now if, as an algebra, $A$ is simple, also $\lambda(A)$ is simple as an algebra. Likewise
as an algebra, $A^{op}$ is simple, and hence also $\rho(A^{op})$.
A classical result \cite[Theorem 4.1.1]{Herstein}
then shows that $\lambda(A)\otimes_K\rho(A^{op})$ is simple.
But now the map
$$A\otimes_KA^{op}\simeq\lambda(A)\otimes_K\rho(A^{op})\lra
\lambda(A)\cdot\rho(A^{op})\subseteq (\End_K^\bullet(C),d_\Hom)$$
is necessarily injective, by the simplicity of $A\otimes_KA^{op}$.


Now, for $n=\dim_K(A)$ we have
$$n^2=\dim_K(A\otimes_KA^{op})=\dim_K(\lambda(A)\cdot\rho(A^{op}))\leq \dim_K(\End_K^\bullet(C,d))=n^2$$
again, and hence
$$A\otimes_KA^{op}\simeq\lambda(A)\cdot\rho(A^{op})=\End_K^\bullet(C,d)$$
as differential graded algebras.
This shows the proposition.
\end{proof}

\medskip


This motivates the following

\begin{Def} \label{equivalentdgalgebras}
Two differential graded $K$-algebras $A$ and $B$,  are called {\em equivalent} 
if there are bounded
complexes $C_1$ and $C_2$ of finite dimensional $K$-modules such that
$$A\otimes_K\End_K^\bullet(C_1)\simeq B\otimes_K\End_K^\bullet(C_2).$$
\end{Def}

We prove now directly the following

\begin{Lemma}\label{Brabelina}
Let $K$ be a field and let $(A,d_A)$ and $(B,d_B)$ be differential graded algebras.
Then $(A\otimes_KB,d_{A\otimes_KB})\simeq  (B\otimes_KA,d_{B\otimes_KA})$ as differential graded algebras.
\end{Lemma}

\begin{proof}
We need to show commutativity of the graded tensor product.
Let $A$ and $B$ be differential graded $K$-algebras.
Then $A\otimes_K B$ has multiplication
$$(a\otimes b)\cdot(a'\otimes b')=(-1)^{|b||a'|}(aa'\otimes bb').$$
The algebra $B\otimes_K A$ has multiplication
$$(b\otimes a)\cdot(b'\otimes a')=(-1)^{|a||b'|}(bb'\otimes aa').$$
Then,
\begin{align*}
A\otimes_K B\stackrel{\alpha}\lra &B\otimes_K A\\
a\otimes b\mapsto&(-1)^{|a||b|}(b\otimes a)
\end{align*}
for any homogeneous elements $a,a'\in A$, $b,b'\in B$
is an algebra isomorphism.
Indeed, for any homogeneous elements $a,a'\in A$ and $b,b'\in B$,
\begin{align*}
\alpha((a\otimes b)\cdot(a'\otimes b'))=&\alpha((-1)^{|b||a'|}(aa'\otimes bb'))\\
=&(-1)^{|aa'||bb'|}\cdot(-1)^{|b||a'|}(bb'\otimes aa')\\
=&(-1)^{(|a|+|a'|)(|b|+|b'|)}\cdot(-1)^{|b||a'|}(bb'\otimes aa')\\
=&(-1)^{|a||b|+|a||b'|+|a'||b'|}(bb'\otimes aa')\\
=&(-1)^{|a||b|+|a'||b'|}(b\otimes a)\cdot(b'\otimes a')\\
=&(-1)^{|a||b|}(b\otimes a)\cdot (-1)^{|a'||b'|}(b'\otimes a')\\
=&\alpha(a\otimes b)\cdot\alpha(a'\otimes b')
\end{align*}
We need to show compatibility with the differential. 
\begin{align*}
\alpha(d(a\otimes b))=&\alpha(da\otimes b+(-1)^{|a|}a\otimes db)\\
=&(-1)^{|b|(|a|+1)}(b\otimes da)+(-1)^{|a|}\cdot(-1)^{|a|(|b|+1)}(db\otimes a)\\
=&(-1)^{|a||b|}(db\otimes a+(-1)^{|b|}b\otimes da)\\
=&d(\alpha(a\otimes b))
\end{align*}
This shows the lemma.
\end{proof}

\subsection{The differential graded Brauer group}

Recall the notion of equivalent simple differential graded
algebras from Definition~\ref{equivalentdgalgebras}.

\begin{Def}  \label{dgbrauerdef}
The {\em dg-Brauer group $\textup{dgBr}(K)$ of
a field $K$} is given by the set of
equivalence classes (in the sense of Definition~\ref{equivalentdgalgebras})
of algebras $A$, which are central simple as $K$-algebras
and which are in addition
differential graded algebras.
\end{Def}

\medskip

Actually, the notion of equivalence from Definition~\ref{equivalentdgalgebras} used
for Definition~\ref{dgbrauerdef}
seems to be a little too strong. We only need to consider complexes $C$ which are
the $K$-module structure of central simple dg $K$-algebras.

\medskip

One might be tempted to define another dg Brauer group by simple dg-algebras 
in the sense of Aldrich and Garcia-Rochas \cite{Tempest-Garcia-Rochas}. This
would then
possibly be linked to the graded Brauer group, which is very well-studied.

\begin{Example} \label{tensorproductofdualnumbers}
Recall Example~\ref{examplesimpledg}. Let $K$ be a field of characteristic $2$ and let
$A=K[X]/X^2$ where $X$ is in degree $-1$, where $d(X)=1$, and $d(1)=0$. Then
$$A\otimes_KA=K[X]/X^2\otimes K[Y]/Y^2=K[X,Y]/(X^2,Y^2)$$
is a differential graded algebra concentrated in degrees $0$, $-1$ and $-2$.
We get that the differential $D$ on the tensor product is
$D(XY)=X-Y$ and $D(X)=1=D(Y)$. Further, $D(1)=0$. This shows that
we have $\ker(D)=K+K(X-Y)$ and hence $\ker(D)=K[Z]/Z^2$ where $Z=X-Y$
is of degree $-1$.
This algebra is not semisimple, whereas $\ker(d)=K$ is semisimple. This shows that
$(A,d)\otimes_K(A,d)$ is not simple as a differential graded algebra in the sense of \cite{Tempest-Garcia-Rochas}, even though $(A,d)$ is simple in the sense of
\cite{Tempest-Garcia-Rochas}. It is therefore impossible to define a Brauer group
for this class of algebras in the same way as it is done classically.
Our observation here confirms our observation from \cite{dgorders} that
the concept of semisimple differential graded algebras in the sense of
\cite{Tempest-Garcia-Rochas} is much less well-behaved as our concept.
\end{Example}

%


\begin{Prop}\label{inverseofdgbr2} 
$\textup{dgBr}(K)$ is an abelian group with group law the graded tensor product and the
inverse element of a class being the equivalence class of the opposite algebra.
\end{Prop}

\begin{proof}
The tensor product is easily seen to be well-defined.
Denote by $[A]$ the equivalence class of a central simple differential graded algebra $A$.
The equivalence class $[K]$ of the $1$-dimensional algebra $K$ concentrated in degree $0$
and differential $0$ is the neutral element of $\dgBr(K)$.
By Proposition~\ref{inverseofA}
for each element $[A]$ of $\dgBr(K)$ we have that
$$[A]\cdot [A^{op}]=[A\otimes_KA^{op}]=[K].$$
Hence all elements of $\dgBr(K)$
have an inverse.
The associativity of the group law is a general property
of the tensor product of differential graded algebras.
By Lemma~\ref{Brabelina} the group law is commutative.

We hence have proved the proposition.
\end{proof}

\section{The main result}

\label{mainresultsection}

\medskip

In \cite{dgorders} we proved a structure theorem for split simple
dg-algebras. In the proof of the theorem we needed a technical hypothesis, namely that there is a
primitive idempotent $e$ of $A$ such that $A\cdot e\not\subseteq A\cdot d(e)$.
In a more general setting this may be false, as is illustrated by the following

\begin{Example} Recall Example~\ref{examplesimpledg}.
Let $K$ be a field, and consider the graded algebra $A=K[X]/X^2$ where $K$ is in
degree $0$ and $X$ is an element of degree $-1$.
Then there is a differential $d$ on $K[X]/X^2$ given by $d(1)=0$ and
$d(X)=1$ such that $(A,d)$ is a differential graded algebra.
Note that here we have $AX\subseteq Ad(X)$. Of course, $A$ is not semisimple and $X$
is not idempotent.
\end{Example}

However, for full matrix algebras over fields, the
hypothesis $A\cdot e\not\subseteq A\cdot d(e)$ in the above mentioned
structure theorem from \cite{dgorders}
is superfluous, as we shall prove now.

Recall the following result by Dascalescu, Ion, Nastasescu, and Rios-Montes
from \cite{DadalescuIonNastacescu}.
Consider the full matrix algebra $\End_K(K^n)$ over a field $K$ and denote by
$e_{i,j}$ the matrix which has coefficient $0$ everywhere
except at position $(i,j)$, where it has coefficient $1$.
There the authors of \cite{DadalescuIonNastacescu} call a group grading on a
full matrix algebra $\Mat_n(K)$  {\em good} if the matrices $e_{i,j}$  are homogeneous
elements of the grading.

\begin{Theorem} \label{DadIonNasta}
\begin{itemize}
\item (\cite[Theorem 1.4]{DadalescuIonNastacescu}).
Let $R$ be the algebra $\Mat_n(K)$ endowed with a $G$-grading such that there is a $G$-graded
$R$-module which is simple as an $R$-module. Then there exists an isomorphism of graded algebras
$R\simeq S$ where $S$ is $\Mat_n(K)$ endowed with a good grading.
\item (\cite[Corollary 1.5]{DadalescuIonNastacescu}). If $G$ is torsion free, then any grading on $\Mat_n(K)$ is isomorphic to a good grading.
\item (\cite[Proposition 2.1]{DadalescuIonNastacescu}).
There is a bijective correspondence between the set of all good $G$-gradings on $\Mat_n(K)$ and the set of maps
$f:\{1,2,\dots,n-1\}\lra G$ such that to a good $G$-grading we associate the map defined by $f(i)=deg(e_{i,i+1})$.
\end{itemize}
\end{Theorem}

We should mention that \cite{DadalescuIonNastacescu} also provide examples of non good gradings
on $\Mat_n(K)$.

The following proposition was proved in \cite{dgorders} under an additional technical assumption.
This assumption is superfluous, as we shall prove now. For the convenience of the reader
we also recall all the details of the proof from \cite{dgorders} in order to have a complete
presentation.

\begin{Prop}\label{progenerator}
Let $K$ be a field
and let $(A,d)$ be a finite dimensional differential graded
$K$-algebra. Suppose that $A$ is a split simple $K$-algebra. Then there is a  bounded
complex $L$ of $K$-modules such that $A\simeq \Hom_K^\bullet(L,L)^{op}$ as differential graded algebras.
Conversely,  $A=\Hom_K^\bullet(L,L)$ is differential graded, finite dimensional
simple as algebra.
\end{Prop}

\begin{proof}
If $L$ is a bounded complex of $K$-vector spaces, then
$\Hom_K^\bullet(L,L)$ is a full matrix ring over $K$, as ungraded algebra, and hence simple
as algebra. Further, as recalled from Section~\ref{generalities} item~(\ref{lastitem}), the algebra
$\Hom_K^\bullet(L,L)$ is a differential graded algebra.

Conversely, let $K$ be a field and let $(A,d)$ be a finite dimensional differential graded
algebra. Suppose that $A$ is a split simple $K$-algebra. By Wedderburn's theorem, $A$ is a
full  matrix algebra over $K$. Let $e$ be a primitive idempotent of $A$.

By Theorem~\ref{DadIonNasta} we may assume that $A$ is $\Z$-graded by a good grading, and
$e_{i,i}$ are all of degree $0$. As $d$ is of degree $1$, we may choose $e=e_{i,i}$
(depending if the
degree $1$ element is upper diagonal or if it is lower diagonal) such that
$Ae\not\subseteq Ad(e)$.

Since $Ae\not\subseteq Ad(e)$,
$$M:=A\cdot e+A\cdot d(e)\textup{ and }N:=A\cdot d(e)$$
are differential graded $(A,d)$-modules. Further $N< M$ and $$L:=M/N\neq 0$$ is a
differential graded $(A,d)$-module.

%
%
%
%
%
%
%

As an $A$-module, we see that $L\simeq Ae$ is a progenerator.
Hence $L$ is a natural differential graded $(A,d)-(\End_K^\bullet(L),d_{\Hom})$
bimodule. Now, for any homogeneous $a\in A$, left multiplication by $a$ gives a
homogeneous element $\varphi(a)\in \End_K^\bullet(L)$. Further, $\varphi$ is additive, sends
$1\in A$ to the identity on $L$, and induces a ring homomorphism
$$\varphi:A\lra  \End_K^\bullet(L)^{op}.$$
Since $L$ is a progenerator, $\varphi$ is injective. Since
$\dim_K(A)=\dim_K(\End_K^\bullet(L))$, we get that $\varphi$ is an isomorphism of algebras.
Now, for any homogeneous $a,b\in A$, we have
$$d(a)b=d(ab)-(-1)^{|a|}ad(b)$$
we get
$$\varphi(d(a))=d\circ\varphi(a)-(-1)^{|a|}\varphi(a)\circ d=d_{\Hom}(\varphi(a))$$
and therefore $\varphi$ is an isomorphism of differential graded algebras.
\end{proof}

{\red\begin{Rem}
After having finished and submitted the manuscript I discovered that D. Orlov
defined our notion of simple differential graded algebra earlier in \cite{Orlov1},
and called it abstractly simple.
Moreover, he proved the statement of Proposition~\ref{progenerator} by completely different means.
His proof uses scheme theoretic arguments. However, he has to assume that the primitive
central idempotents are in degree $0$. Our approach gives that this
can be assumed to be automatically satisfied
using \cite[Corollary 1.5]{DadalescuIonNastacescu}.
\end{Rem}}

\begin{Theorem}\label{derivedbruaergroup}
Let $K$ be a field.
Then the forgetful functor induces an isomorphism
$$\Br(K)\simeq \dgBr(K).$$
\end{Theorem}

\begin{proof}
We obviously have a group homomorphism
$$\Br(K)\stackrel{\iota}\lra \dgBr(K)$$
since any central simple algebra is also a central simple dg-algebra with trivial grading and
$0$ differential.
Further, the map induced by just forgetting the grading and the differential
induces a group homomorphism
$$\dgBr(K)\stackrel{\phi}\lra \Br(K).$$
Of course,
$$\phi\circ\iota=\textup{id}_{\Br(K)}.$$
Hence, $\phi$ is surjective.
Consider $\ker(\phi)$.
By definition $\ker(\phi)$ is formed by equivalence classes of differential graded
central simple algebras $A$ such that $A\otimes_K\Mat_n(K)\simeq \Mat_m(K)$ as
ungraded algebras,
for some $n,m\in\N$. But this shows that $m=n\cdot\dim_K(A)$ and
$$\Mat_n(A)\simeq \End_K(K^{n\cdot\dim_K(A)}).$$
By Proposition~\ref{progenerator} there is a complex $C$ in $C^b(K-mod)$ such that
we have an isomorphism of differential graded algebras
$$\Mat_n(A)\simeq (\End_K^\bullet(C),d_\Hom).$$
Hence $[A]=[K]$ in $\dgBr(K)$. This shows that $\Phi$ is an isomorphism.
Therefore $$\dgBr(K)\simeq \Br(K)$$
and we proved the theorem.
\end{proof}

\begin{Rem}
We emphasize that Theorem~\ref{derivedbruaergroup} shows that if $K$ is a
field, then any central simple differential graded $K$-algebra
$A$ is equivalent to one of the form 
$\End_K^\bullet(C)\otimes_KD$ for some complex $C$ in $C^b(K-mod)$, and some
finite dimensional (ungraded) skew field $D$ with centre $K$.
\end{Rem}

\begin{Rem}
Since $\Br(K)$ is abelian, Theorem~\ref{derivedbruaergroup} shows that
also $\dgBr(K)$ is abelian, without using Lemma~\ref{Brabelina}. However,
Lemma~\ref{Brabelina} is completely elementary, whereas Theorem~\ref{derivedbruaergroup}
is not really. \red{Further, we show in Lemma~\ref{Brabelina} commutativity of the tensor product
in general, and not only up to equivalence in the dg-Brauer group as it follows from 
Theorem~\ref{derivedbruaergroup}.}
\end{Rem}

\begin{Rem}\label{dgbrsecondkindisnotlinkedtoGBr}
Let $(A,d)$ and $(B,d)$ be  simple differential graded algebras in the sense of
\cite{Tempest-Garcia-Rochas}. Since the differential of a tensor product of
dg-algebras $(A,d_A)$ and $(B,d_B)$ is given as $(A\otimes B,d_A\otimes \textup{id}_B\pm \textup{id}_A\otimes d_B)$,
we have the $\ker(d_A)\otimes \ker(d_B)$ is indeed a subspace
of $\ker(d_A\otimes \textup{id}_B\pm \textup{id}_A\otimes d_B)$,
but the kernel contains more in general coming from elements of the form
$x_A\otimes z_Bx_B\pm z_Ax_A\otimes x_B$ for elements $x_A\in\ker(d_A)$ and $x_B\in\ker(d_B)$.
Hence, taking the kernel of the differential will not satisfy that
$(A,d_A)\otimes_K (B,d_B)$ maps to the class of $\ker(d_A)\otimes_K\ker(d_B)$ in $GBr(K)$.
Example~\ref{tensorproductofdualnumbers} is formed in this sense and provides an explicit example
for this phenomenon.
\end{Rem}

\begin{Rem}
By Theorem~\ref{derivedbruaergroup} we have $\Br(K)\stackrel{\iota}{\lra}\dgBr(K)$ is an isomorphism where $\iota([A])=[(A,0)]$.
We may be tempted to consider for a differential graded algebra $(A,d)$ its homology $H(A,d)$.
Since by K\"unneth's formula we get
$$H((A,d)\otimes_K(\End_K^\bullet(C),d_\Hom))\simeq H((A,d)\otimes_KH(\End_K^\bullet(C),d_\Hom))$$
the map 'taking homology' from $\dgBr(K)$ to equivalence classes of graded modules is not really well-defined. Indeed,
$H(\End_K^\bullet(C),d_\Hom)$ is in general not isomorphic to $\End_K(H(C))$.
One would need to consider a broader concept of equivalence for a modified graded Brauer group at least.
Further, even then,
the isomorphism  $\Br(K)\stackrel{\iota}{\lra}\dgBr(K)$ and the fact that $H(A,0)\simeq A$ indicates that
homology is not an interesting map on the differential graded
Brauer group.
\end{Rem}

\end{document}